\documentclass[11pt]{article}
\usepackage{cite}
\usepackage{graphicx}
\usepackage[cmex10]{amsmath}
\usepackage[tight,footnotesize]{subfigure}
\usepackage{amsfonts,latexsym,amssymb,amsmath,amsthm}
\usepackage{url,epstopdf}
\usepackage{tikz}
\usetikzlibrary{arrows,patterns,calc,patterns,decorations.pathmorphing,decorations.markings}
\newcommand{\kron}{\otimes}
\newtheorem{lemma}{\bfseries{Lemma}}
\newtheorem{definition}{\bfseries{Definition}}
\newtheorem{prop}{\bfseries{Proposition}}
\begin{document}
\title{Formation Control of Non-holonomic Vehicles under Time Delayed Communications}
\author{Rudy Cepeda-Gomez\thanks{Institute of Automation, University of Rostock, Rostock, Germany. \texttt{rudy.cepeda-gomez@uni-rostock.de}} and Luis F. Perico\thanks{Faculty of Mechatronic Engineering, Universidad Santo Tomas, Bucaramanga, Colombia. \texttt{luisfernandopericoremolina@mail.ustabuca.edu.co}}}
\maketitle
\begin{abstract}
This document considers the formation control problem for a group of non-holonomic mobile robots under time delayed communications. The agents are assumed to be working under a directed and fixed communication topology. A recently developed formation control technique, based on a consensus protocol for linear systems is revisited. The stability analysis of this control technique, in the domain of the time delays, is performed using the CTCR paradigm under the SDS domain approach. This analysis leads to an exact declaration of the stability boundaries in the domain of the delays. This linear control law is adapted to the nonlinear dynamics of the non-holonomic carts by means of feedback linearization. Simulation results show the applicability of this method to the class of agents considered here.
\end{abstract}

\section{Introduction}
Cooperation among multiple agents is a research topic which has attracted a lot of attention in the past two decades. In order to perform a cooperative task, the members of a team have to reach an agreement on a certain variable. If they want to rendezvous at some point, they must agree upon the location of that point. If they are to travel together, they must select a heading and speed. This gives birth to the \emph{consensus problem}. Since the original agreement protocol of Vickseck et al. \cite{vicsek1995}, expanded by the work of Jadbabaie et al. \cite{Jadbabaie2003}, many authors have worked to extend the theory and applications of the consensus problem. Olfati-Saber and Murray \cite{Olfati-Saber2004} set up a formal framework to treat the agreement problem, highlighting the importance of algebraic graph theory in this topic. They also considered the case in which time delays are present in the communication channels. This pioneering work was limited to agents modeled by single integrator dynamics, and it was soon expanded to consider cases including higher-order dynamics \cite{Ren2005,Ren2008,Lin2008-1,Lin2009} and multiple time delays \cite{Meng2011,IJC2,IJC3}. 

Some researchers have also studied the more general formation generation or flocking problem \cite{Yao2007,Olfati-Saber2006}. Yao et al. \cite{Yao2007} followed an approach based on artificial potentials, whereas Olfati-Saber \cite{Olfati-Saber2006} used a consensus term for velocity agreement and a distance based repulsion term to create the spacing among agents. This methods are heavily non-linear, despite the fact that they are applied to agents with simple single or double integrator dynamics. Following a different approach, Lin et al. \cite{Lin2009} transformed a consensus algorithm into a formation control protocol by including the desired inter-agent distances in the control protocol. 

Cepeda-Gomez and Olgac \cite{IJSS}, also expanded a consensus control law \cite{IJC2,IJC3} to transform it into a formation logic. The last two protocols preserve the linearity of the control logic, and therefore allow the use of well known techniques for the stability analysis with respect to the time delays. All these works, however, have considered agents with very simple dynamics.

Since real life mobile agents have dynamics more complicated than the single or double integrator models commonly used, some authors have expanded the field to include the problem of cooperation of non holonomic agents. The works of Zhai et al. \cite{Zhai2010}, Listman \cite{Listmann2009} and Ghommam \cite{Ghommam2011} present different approaches to solve this problem. However, very few works, besides that of Dong and Farrell \cite{Dong2008}, have treated the consensus or formation control problems for non-holonomic agents under time delayed communications. To study the stability with respect to the delay, they use an approach based on the D-subdivision method, and obtain a delay margin. This result is valid only if all the eigenvalues of the Laplacian matrix generated by the communication topology are real.

In this work, we depart from a consensus protocol affected by two time delays, expand it to include a constant forcing term to create the formation and then deploy it into non-holonomic agents, more specifically, Differential-Drive Mobile Robots, DDMR. A feedback linearization procedure allows us to deploy the original double integrator control protocol to this class of agents, and the Cluster treatment of Characteristic Roots Paradigm \cite{Olgac2005,Fazelinia2007} allows us to find non-conservative stability bounds in the domain of the time delays.

To summarize the contribution of the present paper is the following. Using a feedback linearization model, we adapt a recently developed formation control technique to a set of non-holonomic robots. The formation control technique is based on a consensus protocol and assumes two different time delays in the sensing of positions and velocities. By using the Cluster Treatment of Characteristic Roots methodology, we determine, exactly and exhaustively, the delay combinations that guarantee stability, i.e., convergence to formation in a finite time. To the best of our knowledge, this is the first time such formation problem under time delays is solved in an exactly manner for non-holonomic agents.

The rest of the paper includes the description of the consensus protocol in Section \ref{sec:consensus}, the stability analysis with respect to the delays, using CTCR, in Section \ref{sec:stability}, the formation control algorithm in Section \ref{sec:formation}, the feedback linearization procedure in Section \ref{sec:linearization} and some simulation examples in Section \ref{sec:simulation}. Conclusions and future work directions are presented in Section \ref{sec:conclusions}. Throughout the text boldface notation is used for matrices (uppercase) and vectors (lowercase), whereas italic symbols represent scalar quantities.  
\section{Consensus Protocol}\label{sec:consensus}
We start with the consensus protocol introduced in \cite{IJC3}. It considers a group of $n$ agents with double integrator dynamics of the form $\ddot{\mathbf{x}}_i\left(t\right)=\mathbf{u}_i\left(t\right)$, where $\mathbf{x}_i\left(t\right)\in \mathbb{R}^m$ represents the position vector of the $i$-th agent and $\mathbf{u}_i\left(t\right)\in \mathbb{R}^m$ is the input vector. It is assumed that each one of the position coordinates is independent of the other components of the $\mathbf{x}$ vector, and that the $u_{i,k}$ input affects only the dynamics of the $x_{i,k}$ position, i.e., the dynamics of all the position coordinates are decoupled. Taking advantage of this, in the rest of this section we consider only a scalar position, with the understanding that the results can be easily expanded to the other coordinates for higher-dimensional cases.

We assume the agents are communicating with some of their peers. In the most general case, the communication channels would be unidirectional, and therefore a directed graph is used to model the communication structure. Each vertex represents an agent and each edge a communication channel. The adjacency matrix of this graph, $\mathbf{A}_\Gamma=\left[a_{ik}\right]\in\mathbb{R}^n$, is a binary matrix. The $a_{ik}$ element is 1 if agent $i$ receives information from agent $k$ and 0 if not. The diagonal elements, $a_{ii}$ are taken as 0 too. The $i$-th row sum of $\mathbf{A}_\Gamma$ is called the \emph{in-degree} of the vertex $i$ \cite{Biggs1993}, and it is represented by $\delta_i$. The in-degree is also the number of peers from which the $i$-th member receives information. They are organized in a set called the \emph{informers} of agent $i$, which is denoted by $\mathcal{N}_i$. With the number of informers of each agent, the \emph{in degree matrix} of the communication topology is created. It consists of a diagonal matrix that contains the in-degree of node $i$ in its $i$-th diagonal term: $\boldsymbol{\Delta}=\text{diag}\left(\delta_1,\,\delta_2,\ldots,\delta_n\right)$ 

It is assumed that the communication process is affected by two time delays. The first delay is in the position communication channels and the second one in the velocity information exchange. The different time delays in the velocity and position channels is justified considering that different sensors are used to measure these variables, and therefore different sensing delays will be added on top of the delays inherent to the communication channels. Since we are considering that the members of the group are all alike, assuming the delays are uniform across the network and constant also makes sense. It is also assumed that each agent is able to process its own information fast enough to consider the so called \emph{self delay} negligible. This means, each agent has its current information about position and velocity, but receives only delayed information from its peers. This assumptions are based on the characteristics of the physical platform currently available to perform a practical implementation of the methodology described in this paper.

The control law is defined in this form: each agent takes the \emph{delayed} position information it receives from its informers, computes the average and compares its own \emph{current} position to that average; this error is penalized with a proportional control gain $P\in\mathbb{R}^+$. A similar task is performed with the velocity information received, using a derivative control gain $D\in\mathbb{R}^+$. Mathematically, the control law is represented by:
\begin{equation}
\begin{split}
u_{i}&=P\left(\sum_{k\in\mathcal{N}_i}\frac{x_k\left(t-\tau_1\right)}{\delta_i}-x_i\left(t\right)\right)\\
&+D\left(\sum_{k\in\mathcal{N}_i}\frac{\dot{x}_k\left(t-\tau_2\right)}{\delta_i}-\dot{x}_i\left(t\right)\right),
\end{split}
\label{eq:protocol}
\end{equation}
for $i=1,\,2\,\ldots,\,n$.

Using a concatenation of the positions and velocities of the agents as the state vector $\mathbf{z}=\left[x_1\,\dot{x}_1\,x_2\,\dot{x}_2\,\ldots\,x_n\,\dot{x}_n\right]^T$, the dynamics of the multi-agent system, under the control law  \eqref{eq:protocol}, can be presented in state space as:
\begin{equation}
\begin{split}
\dot{\mathbf{z}}&=\left(\mathbf{I}_{n}\kron\left[\begin{array}{rr}0&1\\-P&-D\end{array}\right]\right)\mathbf{z}\left(t\right)\\
&+\left(\mathbf{C}\kron\left[\begin{array}{rr}0&0\\P&0\end{array}\right]\right)\mathbf{z}\left(t-\tau_1\right)\\
&+\left(\mathbf{C}\kron\left[\begin{array}{rr}0&0\\0&D\end{array}\right]\right)\mathbf{z}\left(t-\tau_2\right)
\end{split}
\label{eq:protocolss}
\end{equation}
where $\mathbf{I}_n$ is the identity matrix of size $n\times n$, $\kron$ represents the Kronecker product operation and the matrix $\mathbf{C}\in\mathbb{R}^{n\times n}$ is the product of the inverse of the in-degree matrix and the adjacency matrix: $\mathbf{C}=\boldsymbol{\Delta}^{-1}\mathbf{A}_\Gamma$.

The complexity of the dynamics in \eqref{eq:protocolss} increases rapidly with the number of agents in the group, quickly becoming intractable. In order to simplify the stability analysis with respect to the time delays, the methodology introduced in \cite{TAC}, extended to the case of directed topologies, is used here. It is based in the following lemma.
\begin{lemma}\label{lemma:factor}
\textbf{Factorization Property}. The characteristic equation of the system \eqref{eq:protocolss} can always be expressed as the product of a $\ell$ second and $m$ fourth order factors as:
\begin{equation}
\begin{aligned}
Q&\left(s,P,D,\tau_1,\tau_2\right)=\prod_{i=1}^{\ell+m}q_i\left(s,P,D,\tau_1,\tau_2,\lambda_i\right)=\\
&\left(\prod_{i=1}^{\ell}\left[s^2+Ds+P-\lambda_i\left(Dse^{-\tau_2\,s}+Pe^{-\tau_1\,s}\right)\right]\right)\\
&\left(\prod_{i=\ell+1}^{\ell+m}\left[s^4+2Ds^3+\left(D^2+2P\right)s^2+2DPs+P^2\right.\right.\\
&-2\Re\left(\lambda_i\right)\left(s^2+Ds+P\right)\left(Dse^{-\tau_2\,s}+Pe^{-\tau_1\,s}\right)+\\
&\left.\left|\lambda_i\right|^2\left(Dse^{-\tau_2\,s}+Pe^{-\tau_1\,s}\right)^2\right]\Bigg)=0
\end{aligned}
\label{eq:factorized}
\end{equation}
where $\lambda_i$ represent the eigenvalues of the $\mathbf{C}$ matrix, which has $\ell$ real eigenvalues, denoted by $i=1,\,2,\ldots,\,\ell$ and $m$ complex conjugate eigenvalue pairs, $\left(\lambda_i,\,\lambda_i^*\right)$, $i=\ell+1,\,\ell+2,\ldots\,\ell+m$. Then $n=\ell+2m$.
\end{lemma}
\emph{Proof} The proof is omitted here due to space constraints. It can be found in \cite{IJC3,IJSS}

As it is shown in \cite{Agaev2005}, if the communication topology is connected and has at least one spanning tree, the matrix $\mathbf{C}$ has 1 as an eigenvalue with multiplicity one. This property introduces a factor in \eqref{eq:factorized} of the form:
\begin{equation}
q_1\left(s,P,D,\tau_1,\tau_2\right)=s^2+Ds+P-\left(Dse^{-\tau_2\,s}+Pe^{-\tau_1\,s}\right)
\label{eq:factor1}
\end{equation}
which always has a root at $s=0$ regardless of the time delays, i.e., it always presents a rigid body mode. This factor is known to dictate the dynamics of a weighted average of the states of the agents, which we call the \emph{weighted centroid} and dictates the value upon which the agents would agree if consensus is reached. Without loss of generality, this eigenvalue is denoted as $\lambda_1$. The other $\ell+m-1$ factors in \eqref{eq:factorized} represent the \emph{disagreement dynamics}: if they are all stable, the consensus is reached. These results are well known in the literature, and are stated here as lemmas for completeness. The proofs of these lemmas can be found in \cite{IJSS}.
\begin{lemma} 
\textbf{Group behavior}. Assume the communication topology has at least one spanning tree. Then, the agents in
the group reach a consensus if and only if the factor \eqref{eq:factor1} is marginally stable and all the remaining factors of \eqref{eq:factorized} are stable.
\end{lemma}

\begin{lemma}
\textbf{Topologies without spanning trees}. If the given communication topology does not have a spanning tree, the control logic described by \eqref{eq:protocol} cannot result in consensus.
\end{lemma}

\section{Stability Analysis Using CTCR}\label{sec:stability}
The results of the previous section show that agents reach consensus if all factors in \eqref{eq:factorized} are stable. Section \ref{sec:formation} shows that this is also a necessary condition to ensure the formation stability. The general form of the factors is:
\begin{equation}
\begin{split}
g_1\left(s\right)+g_2\left(s\right)&e^{-\tau_1\,s}+g_3\left(s\right)e^{-\tau_2\,s}+g_4\left(s\right)e^{-2\tau_1\,s}\\
&+g_5\left(s\right)e^{-2\tau_2\,s}+g_6\left(s\right)e^{-\left(\tau_1+\tau_2\right)\,s}=0
\end{split}
\label{eq:genfactor}
\end{equation}
which is a multiple time delay system with commensuracy (terms like $2\tau_1$) and delay cross-talk (terms like $\tau_1+\tau_2$). In order to find the delay combinations that render the system stable, we deploy the Cluster Treatment of Characteristic Roots paradigm \cite{Olgac2002, Olgac2005}, which we describe next.

It is well known that the characteristic roots of systems like \eqref{eq:genfactor} are continuous functions of the delays, then the stability of this class of systems can only change when there are roots on the imaginary axis. The delay combinations that create imaginary roots can be classified in two classes, defined as follows:
\begin{definition}\label{deff:kernel}
\textbf{Kernel curves} $\wp_0$: The curves that consist of all the points $\left(\tau_1,\tau_2\right)\in\mathbb{R}^{2+}$ exhaustively, which cause an imaginary root $s=j\omega$ and satisfy the constraint $0<\tau_k\omega<2\pi$ are called the \emph{kernel hypercurves}. The points on these curves contain the smallest possible delay values that create the given imaginary root at the frequency $\omega$.
\end{definition}
\begin{definition}
\textbf{Offspring Curves} $\wp$: The curves obtained from the kernel curves by the following pointwise nonlinear transformation:
\begin{equation}
\left\langle\tau_1\pm i\frac{2\pi}{\omega},\,\tau_2\pm k\frac{2\pi}{\omega}\right\rangle,\quad i,\,k=0,1,2,\ldots
\label{eq:nltrans}
\end{equation}
are called the \emph{offspring curves}.
\end{definition}

From a point on the kernel or offspring curves, the toots can move to the left or right half of the complex plane, making the system stable or unstable. This direction of crossing is defined as the root tendency:
\begin{definition}
\textbf{Root Tendency} The root tendency indicates the direction of transition of the imaginary root (to the right or to the left half of the complex plane) as only one of the delays increases by $\epsilon$, $0<\epsilon<1$, while the other delay remains constant: 
\begin{equation}
\left.RT\right|_{s=j\omega}^{\tau_i}=\mathbf{sgn}\left[\mathfrak{Re}\left(\left.\frac{ds}{d\tau_i}\right|_{s=j\omega}\right)\right]
\end{equation}
Root tendencies are –1 for stabilizing and +1 for destabilizing root crossings across the imaginary axis.
\end{definition}
The propositions upon which CTCR is based are the following:
\begin{prop} 
\textbf{Small number of kernel hypercurves}: The number of kernel hypercurves is manageable small: for a LTI-TDS of state dimension $n$, that number is bounded by $n^2$ \cite{Ergenc2007}.
\end{prop}
\begin{prop}
\textbf{Invariant root tendency property}: Take an imaginary characteristic root, $j\omega$, caused by any one of the infinitely many grid points on the kernel and offspring hypercurves in $\left(\tau_1,\tau_2\right)\in\mathbb{R}^{2+}$ domain defined by expression \eqref{eq:nltrans}. The root tendency of these imaginary roots remains invariant so long as the grid points on different \emph{offspring hypercurves} are obtained by keeping one of the delays fixed. That is, the root tendency with respect to the variations of $\tau_i$ is invariant from the kernel to the corresponding offspring as the other delay $\tau_k$ is fixed.
\end{prop}

The proofs of these propositions can be found in \cite{Olgac2002,Ergenc2007,Olgac2005}.

To deploy CTCR, the first step is to find the kernel curves with their corresponding root tendencies. This is followed by the generation of the offspring curves. All these curves partition the domain of the time delays in regions of stability and instability. To find whether a certain region corresponds to stable behavior, we start finding how many unstable roots the non delayed system has. The delays are then increased, and every time a kernel or offspring curve is crossed, the number of unstable roots is increased (if $RT=1$ for that curve) or decreased (if $RT=-1$) by two. Those regions having zero unstable roots, of course, are declared as the stable operating regions.

The question now is how to find the kernel curves. For this, several different methods can be used, like the Kronecker sum \cite{Ergenc2007}, the Rekasius substitution \cite{Olgac2005}, the Spectral Delay Space (SDS) based on Rekasius substitution \cite{Fazelinia2007} or the SDS based on half-angle tangent substitution \cite{IJC2}, which was used in this work. Due to space limitations, we do not include details of the approach we followed here, and we refer the reader to \cite{IJC2,IJC3} where it is explained.
\section{Formation Control}\label{sec:formation}
The consensus protocol in \eqref{eq:protocol} is designed to bring the state of the agents to a common value. If a position state is considered, that would imply a collision. We can transform this behavior into a formation control structure by just adding a constant forcing term to the consensus protocol in \eqref{eq:protocol}. In order to design this term, we must consider the state transformation that decouples the system and factorizes its characteristic equation, as stated in Lemma \ref{lemma:factor} and detailed in \cite{IJC3}. The matrix used in this transformation, $\mathbf{T}$, is the one that takes $\mathbf{C}$ to its Jordan canonical form. Under the transformation $\boldsymbol{\xi}=\left(\mathbf{T}^{-1}\kron\mathbf{I}_2\right)\mathbf{z}$, the dynamics \eqref{eq:protocolss} is decoupled into $\ell$ second order subsystems and $m$ fourth order subsystems. The characteristic equations of these subsystems are the factors of the global characteristic equation as declared in \eqref{eq:factorized}. 

The state space representation of each subsystem, including the extra forcing term we are introducing, has the general form of:
\begin{equation}
\dot{\boldsymbol{\xi}}_i\left(t\right)=\mathbf{A}\boldsymbol{\xi}_i\left(t\right)+\mathbf{B}_{i1}\boldsymbol{\xi}_i\left(t-\tau_1\right)+\mathbf{B}_{i2}\boldsymbol{\xi}_i\left(t-\tau_2\right)+\boldsymbol{\varphi}_i
\label{eq:subsystems}
\end{equation}
where the matrices $\mathbf{A}$, $\mathbf{B}_{1i}$ and $\mathbf{B}_{2i}$ have the following form for subsystems corresponding to real eigenvalues of $\mathbf{C}$:
\begin{equation}
\begin{aligned}
\mathbf{A}=\left[\begin{array}{rr}0&1\\-P&-D\end{array}\right],&\\
\mathbf{B}_{1i}=\lambda_i\left[\begin{array}{cc}0&0\\P&0\end{array}\right],&\ 
\mathbf{B}_{2i}=\lambda_i\left[\begin{array}{cc}0&0\\0&D\end{array}\right]
\end{aligned}
\label{eq:matrices2}
\end{equation}
and the following form for those corresponding to complex conjugate eigenvalue pairs:
\begin{equation}
\begin{aligned}
\mathbf{A}&=\left(\mathbf{I}_2\kron\left[\begin{array}{rr}0&1\\-P&-D\end{array}\right]\right),\\
\mathbf{B}_{1i}&=\left(\mathbf{J}_i\kron\left[\begin{array}{cc}0&0\\P&0\end{array}\right]\right),&\\
\mathbf{B}_{2i}&=\left(\mathbf{J}_i\kron\left[\begin{array}{cc}0&0\\0&D\end{array}\right]\right),&\\
\end{aligned}
\label{eq:matrices4}
\end{equation}
with $\mathbf{J}_i$ being:
\begin{equation}
\mathbf{J}_i=\left[\begin{array}{rr}\mathfrak{Re}\left(\lambda_i\right)&-\mathfrak{Im}\left(\lambda_i\right)\\
\mathfrak{Im}\left(\lambda_i\right)&\mathfrak{Re}\left(\lambda_i\right)\end{array}\right]
\end{equation}
The forcing terms in \eqref{eq:subsystems} have the form $\boldsymbol{\varphi}_i=\left[0\ \phi_i\right]^T$ for second order factors and $\boldsymbol{\varphi}_i=\left[0\ \phi_{1i}\ 0\ \phi_{2i}\right]^T$ for fourth order factors. These forcing terms are applied only to the disagreement dynamics, in order to create the spacing among the agents. The centroid dynamics is left undisturbed, and it dictates where the center of the formation will be located at the end of the exercise.

In order to study the effect of the forcing terms in \eqref{eq:subsystems}, we take the Laplace transform and apply the final value theorem, arriving to
\begin{subequations}
\label{eq:xiss}
\begin{equation}
\xi_{i\infty}=-\frac{\phi_i}{P\left(\lambda_i-1\right)}\left[\begin{array}{c}1\\0\end{array}\right]
\end{equation}
for second order factors and to
\begin{equation}
\xi_{i\infty}=-\frac{\left[\begin{array}{cc}\mathfrak{Re}\left(\lambda_i\right)-1&-\mathfrak{Im}\left(\lambda_i\right)\\0&0\\\mathfrak{Im}\left(\lambda_i\right)&\mathfrak{Re}\left(\lambda_i\right)-1\\0&0\end{array}\right]\left[\begin{array}{c}\phi_{1i}\\\phi_{2i}\end{array}\right]}{P\left(\left|\lambda_i\right|^2+1-2\mathfrak{Re}\left(\lambda_i\right)\right)}
\end{equation}
\end{subequations}
Notice that if we make the forcing terms equal to zero, the steady state disagreement is zero, leading to consensus consensus. When the forcing terms are present, both equations \eqref{eq:xiss} imply constant relative positions among the agents and zero relative velocities.

With the values of $\boldsymbol{\xi}_\infty$, it is possible to obtain the corresponding steady state values for the agents positions, using $\mathbf{z}_\infty=\left(\mathbf{T}\kron\mathbf{I}_2\right)\boldsymbol{\xi}_\infty$. Of course, the formation design process does not start from $\boldsymbol{\xi}_\infty$, but from $\mathbf{z}_\infty$, the desired locations of the agents in the formation. Then, $\boldsymbol{\xi}_\infty$ is found from $\boldsymbol{\xi}_\infty=\left(\mathbf{T}^{-1}\kron\mathbf{I}_2\right)\mathbf{z}_\infty$ and \eqref{eq:xiss} is used to find the value of the forcing terms in the transformed domain, with which the forcing vector $\boldsymbol{\varphi}=\left[0\ 0\ \boldsymbol{\varphi}_2^T\ \cdots\ \boldsymbol{\varphi}_{\ell+m}^T\right]^T$ is created. The actual forcing terms, applied in the non-transformed domain are found using $\mathbf{F}=\left(\mathbf{T}\kron\mathbf{I}_2\right)\boldsymbol{\varphi}$.

With the formation control term added, the dynamics \eqref{eq:protocolss} finally becomes:
\begin{equation}
\begin{split}
\dot{\mathbf{z}}&=\left(\mathbf{I}_{n}\kron\left[\begin{array}{rr}0&1\\-P&-D\end{array}\right]\right)\mathbf{z}\left(t\right)\\
&+\left(\mathbf{C}\kron\left[\begin{array}{rr}0&0\\P&0\end{array}\right]\right)\mathbf{z}\left(t-\tau_1\right)\\
&+\left(\mathbf{C}\kron\left[\begin{array}{rr}0&0\\0&D\end{array}\right]\right)\mathbf{z}\left(t-\tau_2\right)+\mathbf{F}
\end{split}
\label{eq:formation}
\end{equation}
\section{Feedback Linearization}\label{sec:linearization}
The consensus protocol of section \ref{sec:consensus} and the formation control algorithm of section \ref{sec:formation} were devised for second order holonomic agents, in which the motion along different coordinates is independent. This would be the case of some types of omnidirectional mobile robots. In this work, however, we propose the use of non-holonomic wheeled robots, more specifically, of differential drive mobile robots, DDMR, like the \emph{iRobot Create}\cite{createmanual}. In order to apply the protocol of previous sections to this class of robots, some transformation must be undertaken first.
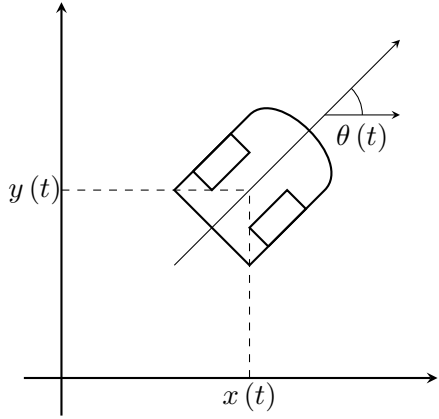
\begin{figure}[t]
\center
\begin{tikzpicture}

\draw [thick] (3.5,2.5) to (2.5,1.5) to (1.5,2.5) to (2.5,3.5) to [out=45,in=45] (3.5,2.5);
\draw [thick] (2.75,1.75) to (3.25,2.25) to (3,2.5) to (2.5,2) to (2.75,1.75);
\draw [thick] (1.75,2.75) to (2.25,3.25) to (2.5,3) to (2,2.5) to (1.75,2.75);
\draw [-stealth] (1.5,1.5) to (4.5,4.5);
\draw [-stealth] (3.5,3.5) to (4.5,3.5);
\draw (4,3.5) arc (0:45:0.5cm);
\draw [thick,-stealth] (-0.5,0) to (5,0);
\draw [thick,-stealth] (0,-0.5) to (0,5);
\node (l1) at (2.5,-0.25) {$x\left(t\right)$};
\node (l2) at (-0.35,2.5) {$y\left(t\right)$};
\node (l3) at (4,3.2) {$\theta\left(t\right)$};
\draw [dashed] (0,2.5) to (2.5,2.5);
\draw [dashed] (2.5,0) to (2.5,2.5);

\end{tikzpicture}

\caption{Variables used to describe the unicycle dynamics.}
\label{fig:uniciclo}
\end{figure}

The kinematics of a DDMR is similar to that of the unicycle. The variables used for its description are the coordinates of the center and the orientation of the wheel, pictured in Fig.~\ref{fig:uniciclo}. Since the wheel is restricted to rolling without slipping, the relationship between the speed, $v$, and the angular rate $\omega$ of the unicycle and the derivatives of its state variables is given by:
\begin{equation}
\left[\begin{array}{c}\dot{x}\\\dot{y}\\\dot{\theta}\end{array}\right]=\left[\begin{array}{cc}\cos\theta&0\\\sin\theta&0\\0&1\end{array}\right]\left[\begin{array}{c}v\\\omega\end{array}\right]
\label{eq:kinem:uni}
\end{equation}

It is well known that \eqref{eq:kinem:uni} does not satisfy Brockett's condition \cite{brockett83}, and therefore there is no continuous, linear time invariant feedback law that can stabilize the system. The consensus protocol \eqref{eq:protocol} cannot be directly applied. However, using a feedback linearization technique, proposed by \cite{Oriolo2002}, we can find an equivalent linear system in which the second derivative of the state equals the control input, as was assumed in Section \ref{sec:consensus}.

To start the feedback linearization of \eqref{eq:kinem:uni}, consider a new state vector $\mathbf{w}=\left[w_1,\,w_2\right]^T$. This vector is selected such that the derivatives of its components are equal to the derivatives of $x$ and $y$:
\begin{equation}
\left[\begin{array}{c}\dot{w}_1\\\dot{w}_2\end{array}\right]=
\left[\begin{array}{c}\dot{x}\\\dot{y}\end{array}\right]=
\left[\begin{array}{c}v\cos\theta\\v\sin\theta\end{array}\right]=
v\left[\begin{array}{c}\cos\theta\\\sin\theta\end{array}\right]
\end{equation}
Taking the second derivative of $\mathbf{w}$ we get:
\begin{equation}
\left[\begin{array}{c}\ddot{w}_1\\\ddot{w}_2\end{array}\right]=
\left[\begin{array}{c}\dot{v}\cos\theta-v\omega\sin\theta\\\dot{v}\sin\theta+v\omega\cos\theta\end{array}\right]
\label{eq:ddz}
\end{equation}
If we make the second derivative of $\mathbf{w}$ equal to the system's input, we obtain the equivalent linear system: 
\begin{equation}
\left[\begin{array}{c}\ddot{w}_1\\\ddot{w}_2\end{array}\right]=
\left[\begin{array}{c}u_1\\u_2\end{array}\right]
\label{eq:syslin}
\end{equation}
By using \eqref{eq:ddz} and \eqref{eq:syslin}, we obtain a relation between the input to the linear system and that of the non-linear system:
\begin{equation}
\left[\begin{array}{c}u_1\\u_2\end{array}\right]=
\left[\begin{array}{cc}\cos\theta&-v\sin\theta\\\sin\theta&v\cos\theta\end{array}\right]\left[\begin{array}{c}\dot{v}\\\omega\end{array}\right]
\label{eq:transform1}
\end{equation}
Inverting the transformation \eqref{eq:transform1} we obtain the input that must be applied to the robot for it to behave as a linear system:
\begin{equation}
\left[\begin{array}{c}\dot{v}\\\omega\end{array}\right]=
\left[\begin{array}{cc}\cos\theta&\sin\theta\\-\frac{1}{v}\sin\theta&\frac{1}{v}\cos\theta\end{array}\right]\left[\begin{array}{c}u_1\\u_2\end{array}\right]
\label{eq:transform2}
\end{equation}
Using \eqref{eq:transform1} and \eqref{eq:transform2}, system \eqref{eq:kinem:uni} is equivalent to \eqref{eq:syslin}, i.e., $x=w_1$, $\dot{x}=\dot{w}_1$, $y=w_2$ y $\dot{y}=\dot{w}_2$.

Care must be taken when using \eqref{eq:transform2}, since it is not defined for $v=0$. This is known to be an structural problem of dynamic extensions of wheeled mobile robots \cite{Oriolo2002}. It should be taken into account when designing control laws for the equivalent linear system.
\section{Simulation Examples}\label{sec:simulation}
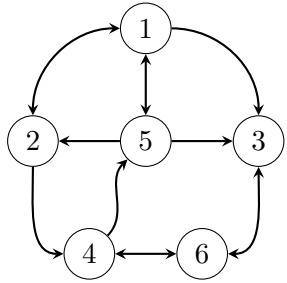
\begin{figure}
\centering
\begin{tikzpicture}
\node (n5) [draw,circle] at(0,0) {5};
\node (n1) [draw,circle] at(0,1.5) {1};
\node (n2) [draw,circle] at(-1.5,0) {2};
\node (n3) [draw,circle] at(1.5,0) {3};
\node (n4) [draw,circle] at(-0.75,-1.5) {4};
\node (n6) [draw,circle] at(0.75,-1.5) {6};
\draw [thick,stealth-stealth] (n1.south) to (n5.north);
\draw [thick,stealth-stealth] (n1.west) to [out=180,in=90] (n2.north);
\draw [thick,-stealth] (n1.east) to [out=0,in=90](n3.north);
\draw [thick,-stealth] (n5.west) to (n2.east);
\draw [thick,-stealth] (n5.east) to (n3.west);
\draw [thick,-stealth] (n2.south) to [out=270,in=180](n4.west);
\draw [thick,stealth-stealth] (n3.south) to [out=270,in=0](n6.east);
\draw [thick,-stealth] (n4.north east) to [out=45,in=225] (n5.south west);
\draw [thick,stealth-stealth] (n4.east) to (n6.west);
\end{tikzpicture}
\caption{Communication topology for the example cases.}
\label{fig:topol}
\end{figure}

In this section we consider 6 agents interacting under the communication topology of Fig.~\ref{fig:topol}. The $\mathbf{C}$ matrix corresponding to this topology has the eigenvalues $\left\{1,\,0.46,\,-0.5,\,-0.64,\,-0.16\pm j0.21\right\}$ which create four second order factors and one fourth order factor, according to \eqref{eq:factorized}. The number of agents was selected because this is the number of actual experimental platforms that are available for a follow-on experimental study. The communication topology was selected such as its $\mathbf{C}$ matrix presented both real and complex eigenvalues without repetitions. For a larger number of agents the analysis methodology remains the same, but the computational load may increase, as larger communication topologies may create more diverse sets of eigenvalues\footnote{A full communication topology, however, is very easy to analyze, as it always creates only two factors, regardless of the number of agents. This fact comes from results from graph theory \cite{Biggs1993} and is shown in detail in \cite{EJC}.}. 

\begin{figure}
\centering
\includegraphics[scale=1.3]{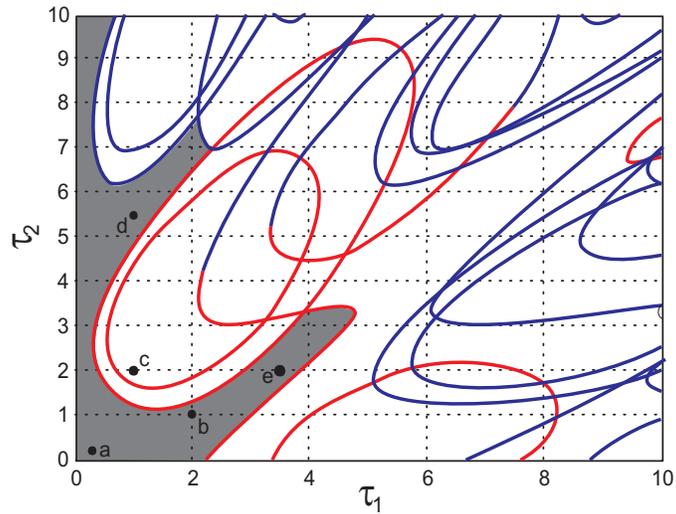}
\caption{Stability chart for six agents interacting under the communication topology of Fig.~\ref{fig:topol} and protocol \eqref{eq:protocol} when $P=1$ and $D=0.5$. The shaded region represents stable delay combinations.}
\label{fig:stability}
\end{figure}

By deploying the CTCR paradigm into the factors, considering control gains $P=1$ and $D=0.5$, we obtain the stability chart of Fig.~\ref{fig:stability} in the domain of the time delays. Each one of the factors creates some kernel, depicted in red, and offspring curves, depicted in blue. At each point on those curves the system has a pair of imaginary characteristic roots. By increasing either of the delays at one of those points the number of unstable characteristic roots changes. The shaded region in that plot corresponds to combinations of delays $\left(\tau_1,\,\tau_2\right)$ for which the number of unstable roots is zero. That is, this is the stable operating region. For any delay combination within the gray area, no factor of the characteristic equation of the system has an unstable root.

We wish to highlight that this chart is exact and exhaustive: no delay combination outside of the gray area results in stable operation, and those in the boundary of the gray area create marginally stable operation, with the agents oscillating around a fixed position.

\begin{figure}
\centering
\includegraphics[scale=0.6]{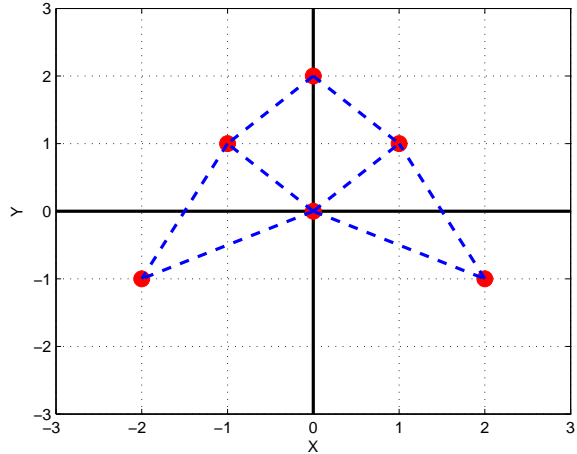}
\caption{Desired formation structure.}
\label{fig:formation}
\end{figure}

\begin{figure}
\centering
\includegraphics[scale=0.8]{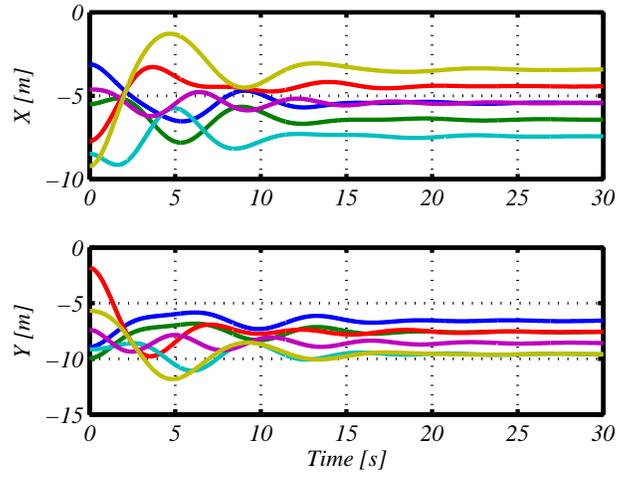}
\caption{Positions of the agents, in $x$ and $y$ direction as a function of time for a delay combination corresponding to point \textbf{a} in Fig.~\ref{fig:stability}.}
\label{fig:xy1}
\end{figure}

\begin{figure}
\centering
\includegraphics[scale=0.8]{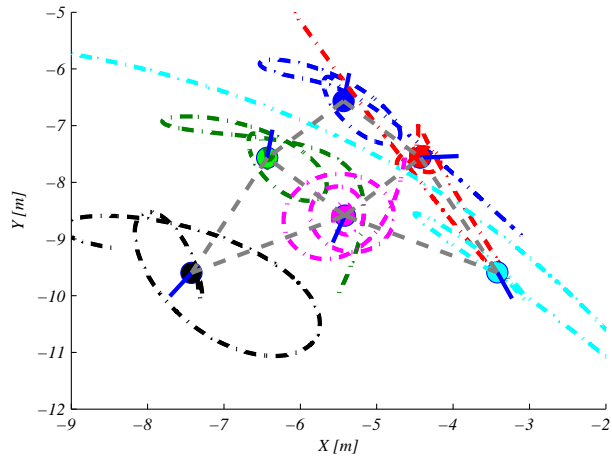}
\caption{Traces of the motion of the agents for a delay combination corresponding to point \textbf{a} in Fig.~\ref{fig:stability}.}
\label{fig:traces}
\end{figure}

To test the formation control algorithm, we propose the formation structure of Fig.~\ref{fig:formation}. For this configuration, the desired $\xi_{\infty}$ values are $[ 0.71,\,1.83,$ $\,1.41,\,-2.08,\,-0.81,\,5.74]$ in the horizontal direction and $[-1.12,\,-0.44,\,1.41,$ $\,1.50,\,-1.40,\,3.72]$ in the vertical direction. In the simulation, the agents start from random positions and orientations, and all have very small ($10^{-3}$) initial speeds to avoid the mentioned division by zero in \eqref{eq:transform2} at $v=0$.
Using the delay combination $\left(\tau_1,\,\tau_2\right)=\left(0.3,\,0.2\right)$, which corresponds to point \textbf{a} in Fig.~\ref{fig:stability}, the agents settle into the formation after about 30 seconds. Figure~\ref{fig:xy1} shows the positions in $x$ and $y$ direction as a function of time, whereas Fig.~\ref{fig:traces} shows how the agents reach the desired formation as final configuration. The small blue lines represent the instantaneous heading of the agents at the moment of the snapshot.

\begin{figure}
\centering
\includegraphics[scale=0.8]{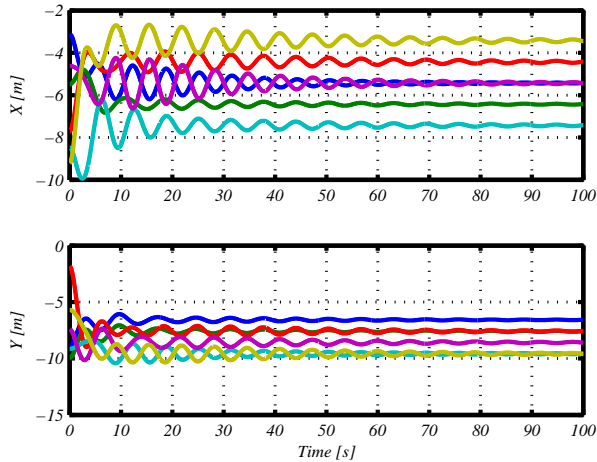}
\caption{Positions of the agents, in $x$ and $y$ direction as a function of time for a delay combination corresponding to point \textbf{b} in Fig.~\ref{fig:stability}.}
\label{fig:xy2}
\end{figure}

An increment of the delays generally worsens the performance of the system. When both delays are increased to $\left(\tau_1,\,\tau_2\right)=\left(2,\,1\right)$, corresponding to point \textbf{b} in Fig.~\ref{fig:stability}, the agents take more than 100 seconds to reach the final state, as can be observed in Fig.~\ref{fig:xy2}.

Increasing the delays, however, is not always bad. An interesting, and rather counter-intuitive, fact that Fig.~\ref{fig:stability} shows is that, in some regions of the delay space, stability can be recovered by increasing one, or both, delays. For example, the delay selection $\left(\tau_1,\,\tau_2\right)=(1,\,2)$, marked by point \textbf{c} in Fig.~\ref{fig:stability}, is in the unstable region. By either increasing $\tau_2$ to 5.5 seconds (point \textbf{d}) or $\tau_1$ to 3.5 seconds (point \textbf{e}), the stability is recovered. This is a very important result, since it provides an extra degree of flexibility to the control designer. Reducing the delays is normally not possible without changing the hardware, but increasing them is simple.

This idea of \emph{delay scheduling} was first discussed in \cite{Olgac2005b,Olgac2007}. It basically consists on a skillful increase of the delays to maintain stability or improve the performance of a system. In the case under study, we would be interested in the delay combination that produces the fastest convergence to the formation. For this, we need to calculate the dominant root of the system, and its associated time constant, for every delay combination within the stable operating zone of Fig.~\ref{fig:stability}. This can be done only with numerical procedures, like the QPmR algorithm \cite{vyhlidal2009}. The stability analysis using the CTCR paradigm is a necessary and crucial first step for delay scheduling, as it is the only methodology that exhaustively declares the stable regions of the system, and therefore provides a finite area in which to search for the smallest time constant.

\begin{figure}
\centering
\includegraphics[scale=0.8]{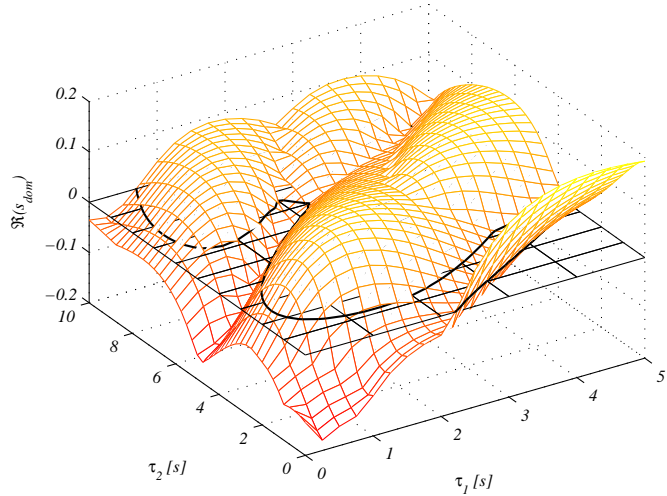}
\caption{Real part of the dominant root of the system as a function of the delays. The black plane corresponds to $\Re \left(s\right)=0$.}
\label{fig:dominant}
\end{figure}

For the system under study, Fig.~\ref{fig:dominant} shows how the real part of the dominant root of the system changes as function of the delays within a domain that includes most part of the stability region. The boundary depicted in Fig.~\ref{fig:stability} can also be seen in Fig.~\ref{fig:dominant} as the intersection of the surface with the plane defined by $\Re\left(s\right)=0$.

\begin{figure}
\centering
\includegraphics[scale=0.8]{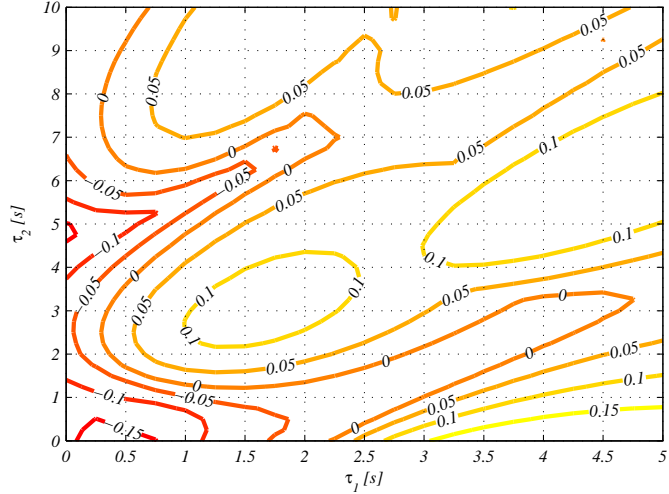}
\caption{Contour plot of the surface in Fig.~\ref{fig:dominant}.}
\label{fig:contour}
\end{figure}

The deepest valleys of the surface in Fig.~\ref{fig:dominant} correspond to the smallest time constants, and therefore to the fastest formation convergence rates. The plot shows that the best delay combinations are those in the vicinity of the origin, and in a small region for $\tau_1\in\left(0,\,1\right)$ and $\tau_1\in\left(4,\,6\right)$. This fact is better observed in Fig.~\ref{fig:contour}, which represents a contour map of the surface in Fig.~\ref{fig:dominant}.

With the information contained in Fig.~\ref{fig:dominant} and Fig.~\ref{fig:contour}, the designer can properly tune the delays to achieve a desired response time for the system. This is the core idea of the delay scheduling.

\section{Conclusions}\label{sec:conclusions}
This paper considered the problem of cooperative formation control for non-holonomic agents in the presence of time delays. We deployed a formation control algorithm based on a linear consensus protocol. We consider two time delays: one in the exchange of position information and the other in the exchange of velocity information. Since the control law was originally devised for linear, holonomic systems, a feedback linearization technique is used to adapt the control logic to the class of agents under study. 

The stability of the formation control depends on the stability of the consensus protocol. In order to find the delay combinations that lead to stable operation, the Cluster Treatment of Characteristic Roots (CTCR) paradigm was used. before using CTCR, a state decomposition that leads to a factorization of the characteristic equation of the system is performed to simplify the problem. With this combination of factorization plus CTCR, it is possible to obtain, in a very efficient manner, an exact and exhaustive stability map with respect to the delays. A simulation shows that the convergence to formation is guaranteed in the regions declared by the proposed methodology.

We also discussed briefly the concept of \emph{delay scheduling} and how the delays can be used to tune the time the agents take to reach the desired formation.

The codes used to create the examples can be requested via email to the authors or can be downloaded from \texttt{https://db.tt/JhOOzQCf}.

Future work includes the practical implementation of these algorithms in a team of mobile robots. For this, the Process Control and Automation Laboratory of the Universidad Santo Tom\'as is working on the development of an experimental platform based on the \emph{iRobot Create} development kit \cite{createmanual}.

\bibliographystyle{IEEEtran}
\bibliography{swarms}
\end{document}